\documentclass{amsart}
\theoremstyle{plain}
\newtheorem{Thm}{Theorem}

\newtheorem{lem}[Thm]{Lemma}

\renewcommand{\theMain}{}

\errorcontextlines=0 \numberwithin{equation}{section}

\newcommand{\rad}{\operatorname{rad}}
\newcommand{\ep}{\epsilon}
\newcommand{\bt}{\bigtriangleup}
\newcommand{\btd}{\bigtriangledown}
\newcommand{\Om}{\Omega}

\begin{document}

\title[Inhomogeneous Boundary Value Problem]
{ Inhomogeneous Boundary Value Problem for Hartree Type Equation}

\author{Li Ma, Pei Cao}

\address{Department of mathematical sciences \\
Tsinghua University \\
Beijing 100084 \\
China} \email{lma@math.tsinghua.edu.cn} \dedicatory{}
\date{Oct 15th, 2009}

\begin{abstract}
In this paper, we settle the problem for time-dependent Hartree
equation with inhomogeneous boundary condition in a bounded
Lipschitz domain in $\mathbb{R}^{N}$. A global existence result is
derived.


\end{abstract}

\thanks{$^*$ The research is partially supported by the National Natural Science
Foundation of China 10631020 and SRFDP 20060003002. }
 \maketitle
 \section{Introduction}

In this paper, we are concerned with the following Hartree type
equation posed on a bounded smooth domain $\Omega$ of
$\mathbb{R}^{N}$
\begin{eqnarray} \label{2}
\left\{
  \begin{array}{ll}
    iu_{t}=\Delta u-\lambda \frac{1}{|x|}*|u|^{2}u, & x\in \Omega, t\in \mathbb{R}\\
    u(x,0)=\phi(x),  \\
    u(x,t)=Q(x,t),&  (x,t)\in \partial \Omega\times \mathbb{R},
  \end{array}
\right.
\end{eqnarray}
where $\phi\in H^{1}(\Omega)$, $\lambda>0$, and $Q\in
C^{3}(\partial\Omega\times(-\infty,+\infty))$ has compact support
and satisfy the compatibility condition on $\partial\Omega$ in the
sense of trace. Here we denote by
$$
\frac{1}{|x|}*|u|^{2}=\int_{\Omega}\frac{1}{|x-y|}|u|^{2}(y)dy.
$$
We now normalize the constant $\lambda=1$.

There is a very large literature on the nonlinear Hartree equation
in the whole space $\mathbb{R}^{N}$, for instance, see
(\cite{CAM}---\cite{MI}). One of the interesting results
concerning existence and uniqueness of $H^1$ global solutions to
Cauchy problem is the paper of Chadam and Glassey \cite{CH}, where
they even show that there exists an unique global solution to the
Cauchy problem of the Hartree system with an external Coulomb
potential. For small data, the scattering of the global solution
is proved in \cite{HA}. Recently,
 $\dot{H}^{1}(\mathbb{R}^{N})$ critical Hartree equation is settled and
 global well-posedness and scattering results are obtained in \cite{MI}.
 Stationary solutions of these equations are of great interest to
 people in the past years. We can solve the remaining problem of E.Lieb that
 the positive solutions to the Hartree equation are minimizers of
 the corresponding functional in \cite{Ma09}.
However, we are aware that very few authors are concerned with the
corresponding problems for Hartree equation in a bounded domain in
$\mathbb{R}^{N}$ with inhomogeneous boundary condition. It is
well-known that inhomogeneous boundary value problems for Hartree
type Schr$\ddot{o}$dinger equations have physical implications.
For instance, in one space dimension, such problem are often
called forced problem when an external force is applied to the
time evolution of systems. Frequently the forcing is put in as a
boundary condition. This is the motivation for us to consider the
problem in this paper.

There are also relatively a few results for non-linear Schrodinger
equations in a bounded domain in $\mathbb{R}^{N}$ with
inhomogeneous boundary conditiion. When space dimension $N>1$,
Strauss and Bu \cite{ST} obtained the global existence of a
$H^{1}$ solution to such Schrodinger equation if the nonlinear
term contributes a positive term to the energy. In \cite{TS},
Tsutsumi obtained interesting result for nonlinear Schrodinger
equation on exterior domains.

In this paper, we establish a global existence result for
(\ref{2}) borrowing an idea from the work \cite{ST}.
 For (\ref{2}), the
conservation of mass and energy do not holds. A direct computation
shows that the growths of the $L^{2}$ norm and energy take the
forms
\begin{eqnarray}\label{5}
\frac{d}{dt}\int_{\Omega}|u|^{2}dx=2Im\int_{\partial\Omega}\overline{u}\frac{\partial
u}{\partial n}dS
\end{eqnarray}
\begin{eqnarray}\label{6}
\frac{d}{dt}E(u)=Re\int_{\partial\Omega}\overline{u}_{t}\frac{\partial
u}{\partial n}dS
\end{eqnarray}
where
\begin{eqnarray}
E(u)=\int_{\Omega}\frac{1}{2}|\nabla
u|^{2}+\frac{1}{4}\frac{1}{|x|}*|u|^{2}|u|^{2}dx
\end{eqnarray}

At first it seems not easy to make use of such identities. Of
course, reducing the problem into homogeneous one does not work
well too. Using more identities derived from the Hartree equation,
we can set up the following existence result.
\begin{Thm}\label{13}
Let $\Omega$ be a bounded smooth domain in $\mathbb{R}^{N}$.
Assume $\phi\in H^{1}(\Omega)$, $Q \in
C^{3}(\partial\Omega\times(-\infty,+\infty))$ have compact support
and  $\phi(x)=Q(x,0)$ on $\partial\Omega$ in the sense of trace.
Then there exists a unique solution $u\in
C_{loc}((-\infty,\infty),H^{1}(\Omega))$ to (\ref{2}).
\end{Thm}

The paper is organized as follows. In the next section, we
introduce some notations and a basic extension result from Sobolev
spaces. Section \ref{sect3} is devoted to apriori estimates of the
solution of (\ref{2}). The local and global existence results are
concerned in Section \ref{sect4}.

\section{Notations and an extension Lemma} \label{sect2}
In this section, we introduce some notations that will be used
throughout the paper.

Assume there exists a smooth real-valued function $\xi$ satisfying
$$
\xi|_{\partial\Omega}=(\xi_{1},\cdots,\xi_{n})=n
$$
where $n$ is the standard outer normal vector for $\partial \Omega$.

Denote by
$$
\partial_{j}=\frac{\partial}{\partial
x_{j}},
$$
$$\eta=\sum_{j}\partial_{j}\xi_{j},
$$
and
$$P=\nabla
u|_{\partial\Omega}.
$$
We also let
$$Re=\text{the real part},\ \ \ Im= \text{the imaginary  part}. $$

We recall a well-known extension result which will play an
important role for apriori estimates later.
\begin{lem}
Assume $\Omega$ is bounded domain in $\mathbb{R}^{N}$ such that
$\partial\Omega$ is Lipschitz, then
 there exists a bounded linear operator
\begin{eqnarray*}
E:H^{1}(\Omega)\longmapsto H^{1}(\mathbb{R}^{N})
\end{eqnarray*}
such that for each $u\in H^{1}(\Omega)$,
\begin{eqnarray*}
E(u)=u  \  \text{a.e. in}\  \ \  \Omega ; \ \
\|Eu\|_{H^{1}(\mathbb{R}^{n})}\leq C\|u\|_{H^{1}(\Omega)}.
\end{eqnarray*}
\end{lem}

For the proof of this result, we refer to the books \cite{Ad} and
\cite{GT}.

\section{apriori Estimates }\label{sect3}

We now employ the idea \cite{ST} to establish apriori estimates
for the smooth solution $u$ of (\ref{2}).
\begin{lem}\label{22}
If $u$ is a smooth solution of (\ref{2}), then $u$ satisfies the
following three identities:
\begin{eqnarray}\label{27}
\frac{d}{dt}(\int_{\Omega}|u|^{2}dx)=2Im\int_{\partial\Omega}\overline{Q}(P\cdot
n)dS,
\end{eqnarray}
\begin{eqnarray}\label{20}
\frac{d}{dt}(\frac{1}{2}\int_{\Omega}|\nabla
u|^{2}dx+\frac{1}{4}\int_{\Omega}\frac{1}{|x|}*|u|^{2}|u|^{2}dx)=Re\int_{\partial\Omega}(P\cdot
n)\overline{Q}_{t}dS,
\end{eqnarray}
and
\begin{eqnarray}\label{18}
&
&\frac{d}{dt}\int_{\Omega}u(\xi\cdot\nabla\overline{u})=2iRe\sum_{j,k}\int_{\Omega}\partial_{k}\xi_{j}\partial_{j}\overline{u}\partial_{k}u
+i\int_{\Omega}\nabla\eta\cdot\nabla\overline{u}u\\
&
&+i\sum_{j}Re\int_{\Omega}\frac{x_{j}}{|x|^{3}}*|u|^{2}\xi_{j}|u|^{2}
+i\int_{\partial\Omega}\frac{1}{|x|}*|u|^{2}|Q|^{2}dS
+i\int_{\partial\Omega}|P|^{2}dS
\nonumber\\
& &-2i\int_{\partial\Omega}|P\cdot
n|^{2}dS+\int_{\partial\Omega}Q\overline{Q}_{t}dS-i\int_{\partial\Omega}\overline{P}\cdot
n\eta QdS.\nonumber
\end{eqnarray}
\end{lem}
\textbf{Proof:} To prove (\ref{27}), we  multiply the equation in
(\ref{2}) by $2\overline{u}$ and then take the imaginary part to
obtain
\begin{eqnarray*}
\frac{\partial}{\partial t}|u|^{2}=2Im\nabla\cdot(\nabla
u\overline{u}).
\end{eqnarray*}
Integrating over $\Omega$ yields to
\begin{eqnarray*}
\frac{d}{dt}\int_{\Omega}|u|^{2}=\int_{\partial\Omega}2Im \nabla
u\cdot n\overline{u}dS=2
Im\int_{\partial\Omega}\overline{Q}(P\cdot n)dS.
\end{eqnarray*}
Similarly,  multiplying the equation in (\ref{2}) by
$2\overline{u}_{t}$ and then taking the real part, we have
\begin{eqnarray*}
2Re\nabla\cdot(\nabla u\overline{u}_{t})-\frac{\partial}{\partial
t}|\nabla u|^{2} + \frac{1}{|x|}*|u|^{2}\frac{\partial}{\partial
t}|u|^{2}=0.
\end{eqnarray*}
Integrate over $\Omega$ to obtain
\begin{eqnarray*}
\frac{d}{dt}\int|\nabla
u|^{2}+\int_{\Omega}\frac{1}{|x|}*|u|^{2}\frac{\partial}{\partial
t}|u|^{2}=2Re\int_{\partial\Omega}P\cdot n\overline{Q}_{t}dS.
\end{eqnarray*}
Since
\begin{eqnarray*}
\int_{\Omega}\frac{1}{|x|}*|u|^{2}\frac{\partial}{\partial
t}|u|^{2}
=\frac{1}{2}\frac{d}{dt}\int_{\Omega}\frac{1}{|x|}*|u|^{2}|u|^{2},
\end{eqnarray*}
thus identity (\ref{20}) follows.

To establish (\ref{18}), noticing
\begin{eqnarray*}
\int_{\Omega}\nabla\cdot(u\xi)\overline{u}_{t}&=&\sum_{j}\int_{\Omega}
\partial_{j}(u\xi_{j})\overline{u}_{t}\\
&=&\sum_{j}\int_{\Omega}\partial_{j}u\xi_{j}\overline{u}_{t}+\int_{\Omega}u\eta\overline{u}_{t}\\
&=&-\int_{\Omega}u\xi\cdot\nabla\overline{u}_{t}+\int_{\partial\Omega}Q\overline{Q}_{t}
dS,
\end{eqnarray*}
we get
\begin{eqnarray}\label{23}
\frac{d}{dt}\int_{\Omega}u(\xi\cdot\nabla\overline{u})&=&\int_{\Omega}u_{t}\xi\cdot\nabla\overline{u}
+\int_{\Omega}u\xi\cdot\nabla\overline{u}_{t}\nonumber\\
&=&\sum_{j}\int_{\Omega}u_{t}\xi_{j}\partial_{j}\overline{u}-\overline{u}_{t}\partial_{j}u\xi_{j}
+\int_{\partial\Omega}Q\overline{Q}_{t}dS-\int_{\Omega}\eta
u\overline{u}_{t}.
\end{eqnarray}
Substitute  the equation in (\ref{2}) into (\ref{23}), we obtain
\begin{eqnarray}
\frac{d}{dt}\int_{\Omega}u(\xi\cdot\nabla\overline{u})&=&
-2iRe\sum_{j}\int_{\Omega}\Delta
u\partial_{j}\overline{u}\xi_{j}+iRe\sum_{j}\int_{\Omega}\frac{1}{|x|}*|u|^{2}\xi_{j}\partial_{j}|u|^{2}\nonumber\\
& &+\int_{\partial\Omega}Q\overline{Q}_{t}dS-\int_{\Omega}\eta
u\overline{u}_{t}.
\end{eqnarray}
Note that
\begin{eqnarray*}
-2iRe\int_{\Omega}\Delta u\partial_{j}\overline{u}\xi_{j}&=&-2iRe\int_{\Omega}\sum_{k}\partial^{2}_{kk}u\partial_{j}\overline{u}\xi_{j}\\
&=&2iRe\sum_{k}\int_{\Omega}\partial_{k}(\partial_{j}\overline{u}\xi_{j})\partial_{k}u-2iRe\sum_{k}\int_{\partial\Omega}\partial_{j}\overline{u}n_{j}\partial_{k}un_{k}dS\\
&=&-iRe\int_{\Omega}\partial_{j}\xi_{j}|\nabla
u|^{2}+iRe\int_{\partial\Omega}n_{j}^{2}|P|^{2}dS\\
&
&+2iRe\sum_{k}\int_{\Omega}\partial_{k}\xi_{j}\partial_{j}\overline{u}\partial_{k}u
-2iRe\sum_{k}\int_{\partial\Omega}\partial_{j}\overline{u}n_{j}\partial_{k}un_{k}dS.
\end{eqnarray*}
Add $j=1,\cdots,N$ to obtain
\begin{eqnarray*}
-2i\sum_{j}Re\int_{\Omega}\Delta
u\partial_{j}\overline{u}\xi_{j}&=&-i\int_{\Omega}\eta|\nabla u|^{2}
+2iRe\sum_{j,k}\int_{\Omega}\partial_{k}\xi_{j}\partial_{j}\overline{u}\partial_{k}u\\
& &+i\int_{\partial\Omega}|P|^{2}dS-2i\int_{\partial\Omega}|P\cdot
n|^{2}dS.
\end{eqnarray*}
Similarly,
\begin{eqnarray*}
i\sum_{j}Re\int_{\Omega}\frac{1}{|x|}*|u|^{2}\xi_{j}\partial_{j}|u|^{2}&=&iRe\sum_{j}\int_{\Omega}\frac{x_{j}}{|x|^{3}}*|u|^{2}\xi_{j}|u|^{2}\\
& &
+i\int_{\partial\Omega}\frac{1}{|x|}*|u|^{2}|Q|^{2}dS-i\int\eta\frac{1}{|x|}*|u|^{2}|u|^{2}.
\end{eqnarray*}
Taking the conjugate of the equation in (\ref{2}), multiplying  by
$\eta \overline{u}$ and integrating over $\Omega$ yields to
\begin{eqnarray*}
\int_{\Omega}\eta u\overline{u}_{t}&=&i\int_{\Omega}\eta
u\Delta\overline{u}-i\int_{\Omega}\eta\frac{1}{|x|}*|u|^{2}|u|^{2}\\
&=&-i\int_{\Omega}\nabla\eta\nabla\overline{u}u-i\int_{\Omega}\eta|\nabla
u|^{2}+i\int_{\partial\Omega}\overline{P}\cdot n\eta
QdS-i\int_{\Omega}\eta\frac{1}{|x|}*|u|^{2}|u|^{2}.
\end{eqnarray*}
Substitute the above results into (\ref{23}), the identity
(\ref{18}) follows.

\section{Local and Global Existence}\label{sect4}

Denote by
\begin{eqnarray*}
f(u)=\frac{1}{|x|}*|u|^{2}=\int_{\Omega}\frac{|u|^{2}(y)}{|x-y|}dy.
\end{eqnarray*}
For Hartree equation (\ref{2}), we need the following Lemma
dealing with the Hartree nonlinearity term.

\begin{lem}
There exists a constant $C>0$ independent of $\|v\|_{H^{1}(\Omega)}$
and $\|w\|_{H^{1}(\Omega)}$ such that
\begin{eqnarray*}
\|f(v)v-f(w)w\|_{H^{1}(\Omega)}\leq
C(\|v\|^{2}_{H^{1}(\Omega)}+\|w\|^{2}_{H^{1}(\Omega)})\|v-w\|_{H^{1}(\Omega)}
\end{eqnarray*}
\end{lem}
\textbf{Proof:} Apply Hardy's inequality in $\mathbb{R}^{N}$
\begin{eqnarray*}
\int_{\mathbb{R}^{N}}\frac{|u|^{2}}{|x|^{2}}dx\leq C\|\nabla
u\|_{L^{2}(\mathbb{R}^{N})}^{2},\forall u\in H^{1}
\end{eqnarray*}
and the extension lemma to obtain
\begin{eqnarray*}
\|f(u)\|_{L^{\infty}(\Omega)}&=
&\sup_{x}\int_{\Omega}\frac{|u|^{2}}{|x-y|}dy\\
&\leq &\int_{\mathbb{R}^{N}}\frac{|Eu|^{2}}{|x-y|}dy \\
&\leq &C\|\nabla Eu\|_{L^{2}(\mathbb{R}^{N})}\|
Eu\|_{L^{2}(\mathbb{R}^{N})}\\
&\leq & C\| Eu\|_{H^{1}(\mathbb{R}^{N})}^{2}\leq
C\|u\|^{2}_{H^{1}(\Omega)}.
\end{eqnarray*}
Similarly ,
\begin{eqnarray*}
\|\nabla f(u)\|_{L^{\infty}(\Omega)}&\leq&
\sup_x\int_{\Omega}\frac{|u(y)|^{2}}{|x-y|^{2}}dy \\
&\leq &
\int_{\mathbb{R}^{N}}\frac{|Eu(y)|^{2}}{|x-y|^{2}}dy\\
&\leq& C\|\nabla Eu\|^{2}_{L^{2}(\mathbb{R}^{N})}\\
&\leq &C\|u\|^{2}_{H^{1}(\Omega)}
\end{eqnarray*}
and
\begin{eqnarray*}
\|f(v)-f(w)\|_{L^{\infty}(\Omega)}&\leq &
\int_{\Omega}\frac{||v|^{2}-|w|^{2}|(y)}{|x-y|}dy\\
&\leq&
C(\|v\|_{L^{2}(\Omega)}+\|w\|_{L^{2}(\Omega)})\|v-w\|_{H^{1}(\Omega)}.
\end{eqnarray*}
For simplicity, we drop the integration domain $\Omega$ hereafter.
For $v,w\in H^{1}(\Omega)$, we have
\begin{eqnarray*}
\|f(v)v-f(w)w\|_{H^{1}}&\leq&C\|f(v)(v-w)+(f(v)-f(w))w\|_{H^{1}}\\
&\leq& C(\|f(v)(v-w)\|_{H^{1}}+\|(f(v)-f(w))w\|_{H^{1}}):=I+II
\end{eqnarray*}
We  now estimate terms I and II by using the above inequalities.
\begin{eqnarray*}
I=\|f(v)(v-w)\|_{H^{1}}&\leq&C(\|f(v)(v-w)\|_{L^{2}}+\|\nabla
f(v)(v-w)\|_{L^{2}}\\
& &+\|f(v)\nabla(v-w)\|_{L^{2}})=C(I_{1}+I_{2}+I_{3}),
\end{eqnarray*}
where
\begin{eqnarray*}
I_{1}=\|f(v)(v-w)\|_{L^{2}}\leq\|f(v)\|_{L^{\infty}}\|v-w\|_{L^{2}}
\leq C\|v\|^{2}_{H^{1}}\|v-w\|_{H^{1}},
\end{eqnarray*}
\begin{eqnarray*}
I_{2}=\|\nabla f(v)(v-w)\|_{L^{2}}\leq\|\nabla
f(v)\|_{L^{\infty}}\|v-w\|_{L^{2}}\leq
C\|v\|^{2}_{H^{1}}\|v-w\|_{H^{1}},
\end{eqnarray*}
and
\begin{eqnarray*}
I_{3}=\|f(v)\nabla(v-w)\|_{L^{2}}\leq\|f(v)\|_{L^{\infty}}\|\nabla(v-w)\|_{L^{2}}\leq
C\|v\|^{2}_{H^{1}}\|v-w\|_{H^{1}}.
\end{eqnarray*}
Similarly,
\begin{eqnarray*}
II=\|(f(v)-f(w))w]\|_{H^{1}}&\leq& C(\|(f(v)-f(w))w\|_{L^{2}}
+\|\nabla(f(v)-f(w))w\|_{L^{2}}\\
& &+\|(f(v)-f(w))\nabla w\|_{L^{2}})=C[I_{4}+I_{5}+I_{6}],
\end{eqnarray*}
where
\begin{eqnarray*}
I_{4}=\|(f(v)-f(w))w\|_{L^{2}}&\leq&
\|(f(v)-f(w))\|_{L^{\infty}}\|w\|_{L^{2}}\\
&\leq& C(\|v\|^{2}_{H^{1}}+\|w\|^{2}_{H^{1}})\|v-w\|_{H^{1}},
\end{eqnarray*}
\begin{eqnarray*}
I_{5}=\|\nabla(f(v)-f(w))w\|_{L^{2}}&\leq &
\|\nabla(f(v)-f(w))\|_{L^{\infty}}\|w\|_{L^{2}}\\
&\leq& C(\|v\|_{H^{1}}+\|
w\|_{H^{1}})\|v-w\|_{H^{1}}\|w\|_{L^{2}}\\
&\leq & C(\|v\|^{2}_{H^{1}}+\|w\|^{2}_{H^{1}})\|v-w\|_{H^{1}},
\end{eqnarray*}
and
\begin{eqnarray*}
I_{6}=\|(f(v)-f(w))\nabla
w\|_{L^{2}}&\leq&\|(f(v)-f(w))\|_{L^{\infty}}\|\nabla
w\|_{L^{2}}\\
&\leq& C(\|v\|^{2}_{H^{1}}+\|w\|^{2}_{H^{1}})\|v-w\|_{H^{1}}.
\end{eqnarray*}
This completes the proof.

The proof of the following lemma is similar to related one given
in \cite{ST}.
\begin{lem}\label{19}
For any $C_{0}>0$, there exists $T_{0}>0$ such that if
$\|\phi\|_{H^{1}}\leq C_{0}$, then there exists a unique solution
$u\in C([0,T_{0}],H^{1}(\Omega))$ which solves (\ref{2}).
\end{lem}
\textbf{Proof:} We choose $\tilde{Q}(x,t)$ to be any $C^{3}$
function on $\overline{\Omega}\times[0,+\infty)$ with compact
support in $x$ such that
$$
\left\{
  \begin{array}{ll}
  \Delta\tilde{Q}=f(Q)Q-iQ_{t}, &  in\   \Omega \\
    \tilde{Q}=Q, & \partial\Omega\times[0,+\infty)
  \end{array}
\right.
$$
Set $v=u-\tilde{Q}$. Then $v$ solves
 \begin{eqnarray}   \label{21}
\left\{
  \begin{array}{ll}
   iv_{t}=\Delta v+\Delta \tilde{Q}-i\tilde{Q}_{t}-f(v+\tilde{Q})(v+\tilde{Q}), & \  in\  \Omega \\
   v(0)=\phi(x)-\tilde{Q}(x,0):=\psi(x),   \\
   v=0 , & \partial\Omega\times[0,+\infty)
  \end{array}
\right.
\end{eqnarray}
Given $\psi\in H^{1}_{0}$, the problem (\ref{21}) can be written as
the integral equation
\begin{eqnarray*}\label{24}
v(t)=e^{-it\Delta}\psi(x)-i\int_{0}^{t}e^{-i(t-s)\Delta}(\Delta
\tilde{Q}-i\tilde{Q}_{t}-f(v+\tilde{Q})(v+\tilde{Q}))ds
:=\mathcal{H}(v)
\end{eqnarray*}
where $e^{-it\Delta}$ is a group of unitary operators on
$H^{1}_{0}(\Omega)$ to itself, $v(t)\in H^{1}_{0}(\Omega)$.\\
Consider the set
\begin{eqnarray*}
E=\{v\in
C([0,T_{0}],H^{1}_{0}(\Omega)):\|v\|_{C([0,T_{0}],H^{1}_{0}(\Omega))}=max_{0\leq
t\leq T_{0}}\|v\|_{H^{1}_{0}(\Omega)}\leq M\}.
\end{eqnarray*}
Define
\begin{eqnarray*}
d(v,w)=\|v-w\|_{C([0,T_{0}],H^{1}_{0}(\Omega))},\  \text{for
every}\ v,w\in E.
\end{eqnarray*}
Thus $(E,d)$ is a Banach space.

 For  $v\in E$,
$\|\psi\|_{H^{1}_{0}}\leq \overline{C}_{0}$ and  each $T_0>0$,
there exist constants $C_{T_{0}}$ and $\tilde{C}_{T_{0}}$ such
that
\begin{eqnarray*}
\|\mathcal{H}(v)\|_{H^{1}_{0}}&\leq
&\|\psi\|_{H^{1}_{0}}+\int_{0}^{t}\|\Delta
\tilde{Q}-i\tilde{Q}_{t}-f(v+\tilde{Q})(v+\tilde{Q})\|_{H^{1}_{0}}dS \\
&\leq &
\|\psi\|_{H^{1}_{0}}+C\int_{0}^{t}(\|v+\tilde{Q}\|^{2}_{H^{1}})\|v+\tilde{Q}\|_{H^{1}}dS +C_{T_{0}}\\
&\leq &\|\psi\|_{H^{1}_{0}}+CT_{0}M^{3}+\tilde{C}_{T_{0}}.
\end{eqnarray*}
Similarly,
\begin{eqnarray*}
\|\mathcal{H}(v)-\mathcal{H}(w)\|_{H^{1}_{0}}&\leq&\|\int_{0}^{t}e^{-i(t-s)\Delta}
(f(v+\tilde{Q})(v+\tilde{Q})-f(w+\tilde{Q})(w+\tilde{Q}))ds\|_{H^{1}_{0}}\\
&\leq&\int_{0}^{t}\|f(v+\tilde{Q})(v+\tilde{Q})-f(w+\tilde{Q})(w+\tilde{Q})\|_{H^{1}_{0}}ds\\
&\leq&\int_{0}^{t}(\|v+\tilde{Q}\|^{2}_{H^{1}}+\|w+\tilde{Q}\|^{2}_{H^{1}})\|v-w\|_{H^{1}_{0}}ds\\
&\leq&2(M+C)T_{0}d(v,w).
\end{eqnarray*}

Take $M=2(\overline{C}_{0}+\tilde{C}_{T_{0}})$ and choose suitable
small $T_{0}>0$ to satisfy $CM^{2}T_{0}<\frac{1}{2}$ and
$2(M+C)T_{0}<\frac{1}{2}$. Then we conclude that $\mathcal{H}$ is
a strict contraction in $(E,d)$. By Banach fixed point theorem,
for any $\overline{C}_{0}>0$, there exists $T_{0}>0$ such that
there is a unique solution $v\in E$ to (\ref{21}).  Thus
$u=v+\tilde{Q}$ is the unique solution in $C([0,T_{0}],H^{1})$ to
(\ref{2}).

Next we use the above results to prove Theorem \ref{13}. In fact,
once establishing the following lemma, Theorem \ref{13} follows.
\begin{lem}\label{28}
Let $T>0$ be fixed. Let $u$ be the solution of (\ref{2}) in the
space $C([0,T],H^{1}(\Omega))$. Then there exists a constant
$C_{T}>0$ such that $\|u\|_{H^{1}}\leq C_{T}$ for all $0\leq t\leq
T$.
\end{lem}
\textbf{Proof:} We remark that we have set up the identities
(\ref{27}), (\ref{20}), (\ref{18}) for the smooth solution of
(\ref{2}) with smooth initial and boundary datum. By the
approximations and a subsequence passage to the limit, we can set
up the existence of the unique solution $u\in C([0,T],H^{1})$ for
general data. Hence we shall always consider smooth data.

At each point on $\partial\Omega$, we can write
\begin{eqnarray*}
|P|^{2}=|P\cdot n|^{2}+|A\cdot\nabla Q|^{2}
\end{eqnarray*}
where $A\cdot\nabla Q$ denotes the tangential component of $P$.
Substituting the above identity into (\ref{18}), we find
\begin{eqnarray*}
& &i\int_{\partial\Omega}|P\cdot
n|^{2}dS+i\int_{\partial\Omega}\overline{P}\cdot n\eta
QdS=-\partial_{t}\int_{\Omega}u(\xi\cdot\nabla\overline{u})
+2iRe\sum_{j,k}\int_{\Omega}\partial_{k}\xi_{j}\overline{u}_{j}u_{k}\\
&
&+i\sum_{j}Re\int_{\Omega}\frac{x_{j}}{|x|^{3}}*|u|^{2}\xi_{j}|u|^{2}
+i\int_{\Omega}\nabla\eta\cdot\nabla\overline{u}u
+\int_{\partial\Omega}Q\overline{Q}_{t}dS\\
& &+i\int_{\partial\Omega}\frac{1}{|x|}*|u|^{2}|Q|^{2}dS
+i\int_{\partial\Omega}|A\cdot\nabla Q|^{2}dS.
\end{eqnarray*}
Note $\xi_{j}$'s are smooth functions on $\Omega$ and $Q\in C^{3}$
with compact support in $x$, then
\begin{eqnarray*}
|\sum_{j}Re\int_{\Omega}\frac{x_{j}}{|x|^{3}}*|u|^{2}\xi_{j}|u|^{2}|
&\leq&\int_{\Omega}\int_{\Omega}\frac{(\sum_{j}(x_{i}-x_{j})^{2})^
{\frac{1}{2}}(\sum_{j}\xi_{j}^{2})^{\frac{1}{2}}}{|x-y|^{3}}|u|^{2}(y)|u|^{2}(x)dydx\\
&\leq&
C\int_{\Omega}\int_{\Omega}\frac{|u|^{2}(y)}{|x-y|^{2}}|u|^{2}(x)dydx\\
&\leq&
C\|\int_{\Omega}\frac{|u|^{2}(y)}{|x-y|^{2}}dy\|_{L^{\infty}(\Omega)}\|u\|_{L^{2}(\Omega)}^{2}\\
&\leq& C\|u\|_{H^{1}(\Omega)}^{4},
\end{eqnarray*}
and
\begin{eqnarray*}
|\int_{\partial\Omega}\frac{1}{|x|}*|u|^{2}|Q|^{2}dS|
&=&|\int_{\partial\Omega}\int_{\Omega}\frac{|u|^{2}(y)}{|x-y|}dy|Q|^{2}(x)dS|\\
&\leq& C\|u\|_{L^{2}(\Omega)}\|u\|_{H^{1}(\Omega)}\\
&\leq& C\|u\|_{H^{1}(\Omega)}^{2}.
\end{eqnarray*}
Combine the above inequalities with the identity (\ref{20}) and
integrate over [0,t] to obtain
\begin{eqnarray}\label{25}
 & &\int_{0}^{t}\int_{\partial\Omega}|P\cdot n|^{2}dSd\tau\leq
|\int_{\Omega}u(\xi\cdot\nabla\overline{u})dx|+|\int_{\Omega}\phi(\xi\cdot\nabla\overline{\phi})dx|\\
& &+C\int_{0}^{t}\int_{\Omega}|\nabla u|^{2}dxd\tau
+C\int_{0}^{t}\|u\|_{H^{1}(\Omega)}^{2}d\tau\nonumber\\
& &+C\int_{0}^{t}\int_{\Omega}|\nabla u||u|dxd\tau
+\int_{0}^{t}\int_{\partial \Omega}|A\cdot\nabla
Q|^{2}dSd\tau\nonumber +C\int_{0}^{t}\|u\|^{4}_{H^{1}(\Omega)}d\tau
\\
& &+\int_{0}^{t}\int_{\partial
\Omega}|Q\overline{Q}_{\tau}|^{2}dSd\tau+C\int_{0}^{t}\int_{\partial
\Omega}|\overline{P}\cdot n Q|dSd\tau.\nonumber
\end{eqnarray}
Denote by
\begin{eqnarray*}
J=(\int_{0}^{t}\int_{\partial\Omega}|P\cdot
n|^{2}dSd\tau)^{\frac{1}{2}}.
\end{eqnarray*}
Since $Q$ is $C^{3}$ with compact support, and $\phi\in
H^{1}(\Omega)$, each term in (\ref{25}) involving $\phi,  Q $ is
 bounded. Therefore
\begin{eqnarray*}
J^{2}\leq c_{0}+c_{1}J+C\int_{0}^{t}\|u\|^{4}_{H^{1}(\Omega)}d\tau
+C\int_{0}^{t}\|u\|^{2}_{H^{1}(\Omega)}d\tau+C\|u\|^{2}_{H^{1}(\Omega)}.
\end{eqnarray*}
Recall the identities (\ref{27}) , (\ref{20})  in Lemma \ref{22}, we
have
\begin{eqnarray*}
\int_{\Omega}|u|^{2}&=&\|u\|^{2}_{L^{2}(\Omega)}\\
&\leq
&\|u(0)\|_{L^{2}(\Omega)}+2Im\int_{0}^{t}|\int_{\partial\Omega}P\cdot
n\overline{Q}dS|d\tau\\
&\leq & c+cJ
\end{eqnarray*}
and
\begin{eqnarray*}
 \int_{\Omega}|\nabla u|^{2}=\|\nabla u\|^{2}_{L^{2}(\Omega)}\leq \int_{\Omega}|\nabla
\phi|^{2}+2\int_{0}^{t}|\int_{\partial\Omega}(P\cdot
n)\overline{Q}_{\tau}dS|d\tau \leq \tilde{c}+\tilde{c}J.
\end{eqnarray*}
So we get by Gronwall's inequality that
\begin{eqnarray*}
J^{2}\leq \tilde c_{1}(T)+\tilde c_{2}(T)J+\tilde
c_{3}(T)\int_{0}^{t}J^{2}(s)dS,
\end{eqnarray*}
which implies that
\begin{eqnarray*}
J\leq C(T)
\end{eqnarray*}
for some uniform constant $C(T)>0$. Therefore $J$ is bounded for
every bounded $T>0$. Thus $\|u\|_{H^{1}(\Omega)}$ is bounded by
$\tilde{C}(T)$. This completes the proof of  Lemma \ref{28}.

For any bounded interval $[0,T]$,  $T=1,2,\cdots$, there exists
unique solution $u\in C([0,T],H^{1}(\Omega))$ to (\ref{2}), and
the case $-\infty\leq t<0$ can be proven in the same way. The
uniqueness follows from the standard argument, one may see
\cite{MA}. Thus we have proved Theorem \ref{13}.

\emph{Acknowledgement}: The authors would like to thank Mr.Jing
Wang for helpful discussion.

\end{document}